\documentclass{article}

\title{On Harmonic Mappings into \\ Weil-Peterssson Completed Teichm\"{u}ller
Spaces}
\author{Sumio Yamada\footnote{Supported in part by NSF DMS0096171/0071862}}

\newtheorem{thm}{Theorem}

\newtheorem{defi}{Definition}

\newenvironment{rmk}{{\bf Remark}}{\\}
\newenvironment{pf}{\smallskip {\bf Proof} \smallskip }{\hfill Q.E.D. \smallskip\\}

\newcommand{\teich}{Teichm\"{u}ller}
\newcommand{\weil}{Weil-Petersson}

\newcommand{\T}{\ensuremath{\cal T}}
\newcommand{\Tbar}{\ensuremath{\overline{\cal T}}}

\newcommand{\g}{\ensuremath{\gamma}}
\newcommand{\del}{\partial}
\newcommand{\map}{\ensuremath{{\rm Map}(\Sigma)}}

\begin{document}

\maketitle

\begin{abstract}
Harmonic mappings into \teich\ spaces appear in the study of 
manifolds which are fibrations whose fibers are Riemann surfaces.
In this article we will study the existence and uniquenesses questions of 
harmonic mappings into \teich\ spaces, as well as
some local and global behavior of the harmonic images induced by 
the \weil\ geometry of \teich\ spaces. 
\end{abstract} 

\section{Introduction}

In~\cite{Y2} the author has studied the geometric structure of the
\weil\ completed \teich\ space of a closed surface of higher genus.  
As explained below it has a stratification
which is respected by the \weil\ isometric action of the 
mapping class group of the surface.  In this article, we will consider 
the behavior of a harmonic mapping into the \weil\ completed \teich\ space
which is equivariant with respect to a representation of the fundamental
group of the domain manifold in the mapping class group.  Note that this 
is a generalization of the classical result of Eells-Sampson~\cite{ES}, where the existence
of a harmonic mapping from a compact domain to a closed target was shown.  The 
analytical difficulty arises because of the lack of the smoothness of the 
target.  However, using the general method developed in~\cite{GS}\cite{KS1}\cite{KS2} and
~\cite{KS3} applied to the specific situation where the target if the \weil\
completed \teich\ space, one can deduce existence theorems.  Once the existence 
is shown, uniqueness questions can also be answered due to the generalized 
negativity of the curvature, as shown in~\cite{Y2} in the sense of Toponogov of the target.  Furthermore
one can say much about how the harmonic image looks within the target,
induced by both the local and the global features of the \weil\ geometry.

Let $\Sigma$ be a closed  topological surface of genus strictly larger than one.
Throughout the paper, we will assume that $\Sigma$ is equipped with a hyperbolic
metric.  In ~\cite{Y2}, the author has studied the geometry of the \weil\ completed 
\teich\ space $\Tbar$ of the surface $\Sigma$ and the \weil\ isometric action of the mapping 
class group $\map$ on $\Tbar$.  It was shown that the completion $\Tbar$ 
is an NPC space (or equivalently CAT(0) space) and and the $\Tbar$ has a
stratification by its boundary components.  To be more precise,
let $\cal S$ be the classes of homotopically nontrivial simple closed curves on the Riemann surface $\Sigma$. Denote by ${\T}_C$ the \teich\ space of (or a product of \teich\ spaces
of) punctured Riemann surface(s) obtained by pinching a collection of {\it mutually
disjoint} simple closed geodesics 
$C = \{ c_i \}$ with $0 \leq i \leq 3g-3$.  Note that $3g-3$ is the upper
limit of the number of mutually disjoint simple closed geodesics on $\Sigma$
of genus $g$.
Then it is shown that 
\[
\Tbar = \cup_{c \in {\cal S}} \overline{{\T}_c} 
\]
where $\T$ is denoted as ${\T}_{\emptyset}$.  It should be noted that $\Tbar$ 
can be also seen as
\[
\Tbar = \cup_{C \subset {\cal S}} {\T}_{C}
\]
since we have the following set theoretic relation
\[
{\T}_{C_1 \cup C_2} \subset \overline{{\T}_{C_1}} \cap \overline{{\T}_{C_2}}
\]
provided $C_1 \cup C_2$ is a subset of $\cal S$ representing  mutually disjoint
simple closed geodesics.  It was shown in ~\cite{Y2} that for each $C \subset \cal S$
of mutually disjoint simple closed geodesics,
${\T}_C$ is totally geodesic in the sense that given a pair of points in 
${\T}_C$, the unique distance-realizing \weil\ geodesic connecting
them lies entirely in ${\T}_C$.  Also any open \weil\ geodesic segment in $\T$
lies entirely in ${\T}_C$ for one unique $C$ in $\cal S$, which can be 
stated as  any two 
boundary components of the same dimension which meet do so transversely, 
in fact it is a normal intersection. 

Let ${\cal D}$ be the group
of smooth diffeomorphisms of $\Sigma$ , and ${\cal D}_0$ the group of diffeomorphisms
homotopic to the identity map from a fixed Riemann surface $\tilde{\Sigma}$
This gives a {\it marking} to all the points in ${\cal M}_{-1}$.
The discrete group ${\cal D}/{\cal D}_0$ is called the mapping class group,
or the \teich\ modular group. which we will denote by $\map$.  It is well known 
that $\map$ is generated by Dehn twists. Also it should be noted that  the mapping class group 
(or equivalently \teich\ modular group)
is contained in  the \weil\ isometry group of the \teich\ space.  
It was shown in ~\cite{Y2} that
the isometric action of $\map$ on $\T$ extends to $\Tbar$, due to the fact that 
the displacement distance $d(x. \g x)$ by a Dehn twist \g\ around a simple closed curve 
is zero when $x$ is a point in the boundary component ${\T}_c$.  

Therefore we have this NPC space $\Tbar$ with the mapping class group $\map$ acting isometrically.  

In this paper, we will address the
existence and uniqueness  
of  harmonic maps ({\bf Theorem 1, 2} and {\bf  3}) into the \weil\ completed \teich\ space
as well as how the images of harmonic maps 
respect the stratification structure of $\Tbar$ ({\bf Theorem 4} and {\bf 5}). Those results
have applications in studying the geometry of Lefschetz pencil/fibration as demonstrated 
below.

\section{Statements of Theorems}

We will first define a functional defined on the \weil\ completion $\Tbar$
of the \teich\ space $\T$.  We denote by $d(x, y)$ the \weil\ distance 
between the points $x$ and $y$.  

\begin{defi}
Suppose $\Gamma$ is a finitely generated subgroup of $\map$, with
$\{ {\g}_i \} 1\leq i \leq l$ its generators.  The we define a functional
$\delta$ on $\Tbar$ by
\[
\delta(x) = \max_{1 \leq i \leq l} d(x, {\g}_ix).
\]
\end{defi}

Note here that $\delta: \Tbar \rightarrow {\bf R} \cup + \infty$
is a convex functional, since each $d(x, {\g}_ix)$ is convex
on $\Tbar$ due to the NPC curvature condition, and since
the maximum of finitely many convex functionals is again convex.  

\begin{defi}
Given a subgroup $\Gamma$ of $\map$, the isometric action of $\Gamma$
on $\Tbar$ is said to be {\it proper} if the sublevel set
\[
S(M) = \{ x \in \Tbar : \delta(x) < M < + \infty \} 
\]
is bounded in $\Tbar$.
\end{defi}

The statement of the existence theorem now follows.\\
\\
\noindent {\bf Theorem (Existence Theorem 1)}
{\it Suppose $M$ is a compact manifold without boundary.  Suppose that a
representation $\rho : \pi_1 (M) \rightarrow {\rm Isom}( {\Tbar} )$ is such that 
the image $\rho[\pi_1(M)] $ in ${\rm Isom}( {\Tbar} )$ induces a proper action
on $\Tbar$.  Then there exists an energy minimizing harmonic map 
$u : \tilde{M} \rightarrow \Tbar$ which is $\rho$-equivariant.
($\tilde{M}$ is the universal covering space of $M$.) Moreover the 
map $u$ is Lipschitz continuous. 
}
\\

When $M$ is not compact but complete, under mild additional conditions which 
often are met for applications, we still have the existence theorem.\\
\\
\noindent {\bf Theorem (Existence Theorem 2)}
{\it Suppose that $M$ is a noncompact, complete,  that there exists
a finite energy $\rho$-equivariant map $\phi : \tilde{M} \rightarrow
\Tbar$, and that $\rho$ induces a proper action on $\Tbar$.  Then
there exists an energy minimizing harmonic map $u : \tilde{M} \rightarrow \Tbar$  
which is locally Lipschitz continuous.  
}
\\
As for uniqueness of the map thus obtained, the answer is affirmative.\\
\\
\noindent {\bf Theorem (Uniqueness of Harmonic Map)}
{\it The  harmonic map is unique
within the class of finite energy maps which are $\rho$-equivariant, provided
that the image of the map is not contained in a geodesic.
}\\

In order to state a further consequence of the existence/uniqueness theorems  we first need to recall 
the classification theorem~\cite{Th} of elements of the mapping class group $\map$ first and 
then that of subgroups of  $\map$~\cite{MP}.  

An element of $\map$ is classified as one of the following three types:\\
1) it is of finite order, also called periodic or elliptic;\\
2) it is reducible if it leaves a tubular neighborhood of 
a collection $C$ of closed geodesics $c_1,... c_n$ invariant;\\
3) it is pseudo-Anosov (also called irreducible) if there is $r>1$ and transverse measured
foliations $F_+, F_{-}$ such that $\g (F_{+}) = r F_{+}$ and 
$\g (F_{-}) = r^{-1} F_{-}$.  In this case the fixed point set of
$\g$ action in ${\cal PMF}( \Sigma)$ (the Thurston boundary of $\T$) 
is precisely $F_+, F_{-}$.\\

As for classification of subgroups, McCarthy and Papadapoulos~\cite{MP} have shown
that  the subgroups of $\map$ is classified into four classes:\\
1) subgroup containing a pair of independent pseudo-Anosov elements
(called {\it sufficiently large subgroups};\\
2) subgroups fixing the pair $\{ F_+ (\g ), F_{-} (\g ) \}$ of fixed points
in ${\cal PMF} (\Sigma)$ for a certain pseudo-Anosov element $\g \in \map$ (such groups
are virtually cyclic);\\
3) finite subgroups;\\
4) infinite subgroups leaving invariant a finite, nonempty, system of
disjoint, nonperipheral, simple closed curves on $\Sigma$
(such subgroups are called reducible.)

We are now ready to state the following result.\\
\\
\noindent {\bf Theorem (Stratification Property of Harmonic Images)}
{\it Suppose that $\rho$ induces a proper action on $\Tbar$ and that there
exists a unique energy minimizing harmonic map $u: \tilde{M} \rightarrow \Tbar$.
then the subgroup $\Gamma = \rho (\pi_1 (M))$ is reducible by a collection $C$
of mutually disjoint closed geodesics, if and only if the image of the energy minimizing map
$u$ lies entirely in the boundary component $\overline{{\T}_C}$. }

\begin{rmk}
Note here that when $C$ as a subset of $\cal S$ is the empty set, this is just the
statement of the existence theorems above.  What this result says is that the reducibility
of the subgroup $\Gamma$ of $\map$ is equivalent to the reducibility of the $\Gamma$-invariant 
harmonic map in accordance with  the stratification structure of the NPC space $\Tbar$.
\end{rmk}

We will discuss an application.
A K\"{a}hler four manifold is said to have a structure of a holomorphic
Lefschetz fibration if the following descriptions are met.
There exists a holomorphic map
$\Pi : M \rightarrow B$ where $B$ is a closed Riemann surface.
the map $\Pi$ has finitely many critical points $N_i, \  i=1,..., n$ 
in disjoint fibers $F_i = \Pi^{-1} (P_i), \  i+1, ..., n$. 
The neighborhood of each critical point $N_i$ 
can be described locally by complex coordinates  $z, w$ on $M$ and 
$t$ on $B$ such that $\Pi: (z, w) \rightarrow t(z, w) = zw$ where 
$N_i = (0, 0)$ and $P_i = 0$.  In other words, 
the fiber $F_i$ has a node $N_i$.  
Away from those disjoint fibers, each  fiber of the 
map $\Pi$ is a Riemann surface of varying conformal structures of a fixed genus
$g$. 

The picture above can be transcribed as saying that there exists a 
$\rho$-equivariant holomorphic map $\Pi :B \backslash \{ P_i \}
\rightarrow \T$, where $\rho: \pi_1 (B \backslash \{ P_i \}) 
\rightarrow \map $ is the monodromy representation of the fibration.

The stratification theorem above has an immediate application, which had been first proven by Shiga~\cite{Sh} by a 
different method.\\
\\
\noindent {\bf Corollary}
{\it Given a holomorphic Lefschetz fibration of higher genus, its monodromy representation is sufficiently
large.}

\begin{pf}
It is well known that a holomorphic map between two K\"{a}hler 
manifolds is energy-minimizing~\cite{Sa}.  Hence the map $\Pi
: B \backslash \{ p_i \} \rightarrow \Tbar$ 
is the unique $\rho$-equivariant harmonic map whose existence
and uniqueness have been so far established.  To see that $\Gamma$
is sufficiently large, note if it weren't, then we have three other possible
cases.  The first being when $\Gamma$ is a finite group can be
excluded since each local monodromy is of infinite order.
The second being the case that $\Gamma$ is virtually cyclic, fixing
a pair of points in the Thurston boundary.  Then the image of
the $\rho$-equivariant harmonic map $u : \tilde{B \backslash \{ p_i \}} 
\rightarrow \Tbar$
is a $\Gamma$ invariant \weil\ geodesic, which lies entirely in 
the interior \teich\ space $\T$.  The projection of the invariant
geodesic down to the moduli space is a loop located 
away from any of the divisors, which in turn says that there is no
sequence of points $\{ q_j \}$ in $B \backslash \{p_i\}$ 
over which a cycle on 
the Riemann surface represented by $u(p_i)$ is vanishing
(or equivalently a neck is pinching), which contradicts with the
fact that $M^4$ has singular fibers/vanishing cycles.  

Lastly the third possible case to be excluded is 
when $\Gamma$ can be reduced by a collection
of $C$ mutually disjoint simple closed curves.  
Then the harmonic image of $\Pi$ is entirely contained
in $\overline{{\T}_C}$, which implies that every fiber is a Riemann
surface with nodes where the nodes are obtained by pinching each
simple closed curves in $C$.  
This certainly is not the case when $M$ is a Lefschetz fibration. 
 
\end{pf}

Note that the previous theorem says that given the equivariance condition, the
harmonic image is always found within the subset of $\Tbar$ where the representation
of the fundamental group is no longer reducible.  Hence in a sense it gives an upper
bound on the size of the set an energy minimizing map maps into.  The following theorem, on the other
hand, gives a lower bound.  Namely any harmonic image is contained within a 
copy of \teich\ space in $\Tbar$ (and it does not ``spill over'' onto the boundary.)

 First we introduce a class of maps we will consider.

\begin{defi}
A map into the \weil\ completed \teich\ space is said to be geometrically
finite if for each point $p$ in the domain, there exists a neighborhood
$\Omega$ of $p$ such that there exists a set $\{ c_i \}_{i=1}^{3g-3}$ of mutually disjoint set of simple closed geodesics
on the surface such that the image of the map restricted to $\Omega$ is contained in a bounded set with respect
to the Fenchel-Nielsen coordinates $(l_i, {\theta}_i)$ of $\{ c_i \}$.
\end{defi}

\begin{rmk} 
Note that given a map of finite energy, it is clear that the image is bounded
in the hyperbolic length parameters $l_i$,$( 1 \geq i \geq 3g-3)$ of 
the Fenchel-Nielsen coordinate.
However, the finite energy condition does not imply the boundedness of 
the image in the Fenchel-Neilsen twist parameters ${\theta}_i$, since
the \weil\ metric tensor degenerates as a node develops~\cite{Ma}.
\end{rmk}

\noindent {\bf Theorem (Unique Continuation Property of Harmonic Images)}
{\it Suppose that a Riemannian region $\Omega$ is mapped into ${\Tbar}$ 
by a geometrically finite energy-minimizing map $u$, and that there is a point 
$p$ in $\Omega$ mapped into a \teich\ space ${\T}_C \subset {\T}$ for some collection of
mutually disjoint simple closed curves, where C possibly is the empty set.  Then  the entire
image $u[\Omega]$ lies in ${\T}_C$.}

\begin{rmk}
A Riemannian region is a simply connected open susbet of ${\bf R}^n$ with a
Riemannian metric.
\end{rmk}

\begin{rmk}
It should be noted that the statement of the theorem is a higher dimensional 
generalization of the one dimensional result ({\bf Theorem 3}) in~\cite{Y2}; {\em Every
open \weil\ geodesic segment in $\Tbar$ is entirely contained in a
single copy of \teich\ space ${\T}_C$ for some $C$ (possibly the empty set.)}
\end{rmk}
\\
\noindent {\bf Corollary (Harmonic Map induced by Lefschetz Fibration induced Monodromy)}
{\it
Suppose that $M^4$ is a closed four manifold which fibers holomophically
over $S^2$ with $N$ singular fibers over $\{ q_i \}_{i = 1}^n$ in $S^2$.
in the sense of Lefschetz.  Let $\rho: \pi_1 (S^2 \backslash \{ q_i \})
\rightarrow \map$ to be the monodromy representation of the holomorphic 
Lefschetz fibration.  Suppose now that for a set of $N$ distinct points
$\{ p_i \}_{i=1}^N$ on $S^2$ there exists a $\rho$-equivariant energy minimizing 
map $u : \widetilde{S^2 \backslash \{p_i \}} \rightarrow \Tbar$.  Then the image
of $u$ lies in $\T$.   
}

\begin{pf}[of the Corollary]  According to the theorem above, it suffices to prove 
that $u$ restricted to a fundamental region $F$ in its universal covering
space $\widetilde{S^2 \backslash \{p_i \}}$ is geometrically finite. For if the map
is geometrically finite, and if a point $p$ gets mapped into ${\T}_C$ for some nonempty
$C$, then the theorem above says that the entire image of $u$ is contained in
${\T}_C$, which in turn is only possible when the subgroup of $\map$ representing
the monodromy of the holomorphic fibrations is reducible.  However, the
representation $\rho$ realizes a holomorphic Lefschetz fibration $M^4$, which cannot
be reducible by the previous theorem above. a contradiction.  Therefore $C$ above
has to be empty, and hence the image of $u$ is entirely contained in $\T$.

Looijenga~\cite{Lo} showed that there exists a finite index
normal subgroup $\Gamma$ of $\map$ such that ${\T} / \Gamma$ is a smooth
manifold, and its Deligne-Mumford compactification ${\Tbar} / \Gamma$
is a smooth variety,
where its quotient space by the finite group
${\map} / \Gamma$ is the Deligne-Mumford compactification
of the moduli space $M_g$.  Let $K$ be the kernel of the $\rho$-induced
representation $ {\pi}_1 (S^2 \backslash \{p_i\}) \rightarrow {\map} / \Gamma$.
Let $M$ be the quotient space $\widetilde{S^2 \backslash \{p_i\}}  / K$, which is
a finite cover of $S^2 \backslash \{p_i\}$.  Then $u$ descends to a map
$\overline{u} : M \rightarrow {\Tbar} / \Gamma$, where the target is a smooth
variety.  By the compactness of ${\Tbar} / \Gamma$, it follows that 
$\overline{u}$
is geometrically finite.  Since $M$ is a finite cover of $S^2 
\backslash \{p_i\}$
it follows that $u$ is geometrically finite on a fundamental domain $F \subset 
\widetilde{S^2 \backslash \{p_i\}}$.
\end{pf}

\section{Proofs of the Existence Theorems}

We will present the proofs of the two existence theorems here.  The general
theory of harmonic mappings into NPC spaces was developed in a series of 
papers~\cite{KS1} \cite{KS2} \cite{KS3} by Korevaar and Schoen,
and it is their argument that is adapted here in our context.  

\begin{thm}[Existence Theorem 1]
Suppose $M$ is a compact manifold without boundary.  Suppose that a
representation $\rho : \pi_1 (M) \rightarrow {\rm Isom}( {\Tbar} )$ is such that 
the image $\rho[\pi_1(M)] $ in ${\rm Isom}( {\Tbar} )$ induces a proper action
on $\Tbar$.  Then there exists an energy minimizing harmonic map 
$u : \tilde{M} \rightarrow \Tbar$ which is $\rho$-equivariant.
($\tilde{M}$ is the universal covering space of $M$.) Moreover the 
map $u$ is Lipschitz continuous.
\end{thm}

\begin{pf}
  We will first define the space $Y = L^2_{\rho} (\tilde{M}, \Tbar)$ of
$\rho$-equivariant locally $L^2$ maps.  The space $Y$ has a natural distance 
function $d_2$ given by 
\[
d_2(u, v) = \Big( \int_{\tilde{M}/\Gamma} d_{wp} (u(x), v(x)) d\mu \Big)^{1/2}
\]
where we note that $x \rightarrow d_{wp} (u(x), v(x))$ is $\rho$-equivariant
on $\tilde{M}$, and hence defines a function $M = \tilde{M}/\Gamma$.  Recall that
the fact that $\Tbar$ is an NPC space implies that for almost any $x$, we have
\[
d^2_{wp} (w(x), u_{1/2}(x)) \leq \frac{1}{2} d^2_{wp} (w(x), u_0(x))^2 + 
\frac{1}{2} d^2_{wp} (w(x), u_1(x))^2 - \frac{1}{4} d^2_{wp} (u_0(x), u_1(x))^2
\]
where $u_{1/2}(x)$ is the mid-point on the \weil\ geodesic connecting $u_0(x)$ and
$u_1(x)$ in $\Tbar$.  Integrating the inequality over $\tilde{M}/\Gamma$, we get
\[
d^2_2 (w, u_{1/2}) \leq \frac{1}{2} d^2_2 (w, u_0)^2 + 
\frac{1}{2} d^2_2 (w, u_1)^2 - \frac{1}{4} d^2_2 (u_0, u_1)^2
\] 
which implies that $\Tbar$ is NPC.  Note that the geodesic connecting 
a pair of points in $Y$ is given by taking the geodesic homotopy interpolating
the two points/maps. 

By~\cite{KS1} there exists a sequence $\{ u_i \}$ of $\rho$-equivariant maps with
$E(u_i) \geq E(u_{i+1})$, and 
\[
\lim_{i \rightarrow \infty} E(u_i) = E_0
\]
where $E_0 = \inf \{ E(u): \mbox{$u$ $\rho$-equivariant} \}$.
Then it was shown~\cite{KS1} that the sequence $\{ u_i \}$ is locally uniformly Lipschitz (
a generalized DiGiorgi-Nash estimate.)  Let $C_i$ be the closed convex hull in $Y$ of
the tail end of the sequence.  The sequence $\{ C_i \}$ is a nested sequence of closed
convex sets of $Y$.  Since the energy is convex with respect to $d_2$ (due to the fact 
that $\Tbar$ is NPC) on the subset of $Y$ consisting of finite energy maps, and is 
lower semi-continuous with respect to $d_2$ convergence, we have 
\[
\sup_{u \in C_i} E(u) \leq E(u_i) 
\]
for all $i=1, 2, ...$. We have also used here that $\{ E(u_i) \}$ is a decreasing sequence.
  If there is a common point in all of the $C_i$, then $u$ is 
a harmonic equivariant map.  Since the intersection of a nested decreasing sequence
of nonempty closed bounded convex sets in an NPC space is nonempty (\cite{KS2}),
we now want to show that $C_1$ is bounded, in order to show the existence.  

Denote by $F$ the fundamental region in $\tilde{M}$ for the action of ${\pi}_1 (M)$ on $\tilde{M}$.
Let $P_0$ be a chosen point,and let $u_0$ be the equivariant  map such that 
$u_0(x) = P_0$ for $x \in F$.  Because $\overline{F}$ is compact and the $\{ u_i \}$ are
 locally Lipschitz uniformly in $i$, we have 
\begin{eqnarray}
d_{wp} (u_i(\gamma_{\alpha} x), u_i(x)) & = & d_{wp} ( \rho (\gamma_{\alpha}) u_i(x), u_i(x))\\
 & \leq  & L d_{\tilde{M}}(\gamma_{\alpha}x, x) \\
 & < & C.
\end{eqnarray}
for  $x \in F$,  $\alpha = 1, 2, ...n$, and some $0< C < \infty$, where $L$ is the Lipschitz
bound of the maps $u_i$. 
Recall the functional $\delta : \Tbar \rightarrow {\bf R}$ 
was defined to be $\delta(x) = \max_{\alpha = 1, ... l} d_{wp} (\gamma_{\alpha} x, x)$.
We have just observed $\delta (u_i (x)) < C$ for $x \in F$ and for all $i$.
Since the action is proper, we know that the set $u_i [F]$ is contained in the
bounded set $\{ x \in \Tbar : \delta (x) < C \}$.  It now follows that $d_2(u_i, u_0) < C'$
for all $i$ and therefore 
\[
d_2 (u, u_0) < C'
\]
for all $u \in C_1$.

\end{pf}
\begin{thm}
Suppose that $M$ is a noncompact, complete,  that there exists
a finite energy $\rho$-equivariant map $\phi : \tilde{M} \rightarrow
\Tbar$, and that $\rho$ induces a proper action on $\Tbar$.  Then
there exists an energy minimizing harmonic map $u : \tilde{M} \rightarrow \Tbar$  
which is locally Lipschitz continuous.
\end{thm}

\begin{pf} We will quote the following theorem from~\cite{KS1}.

\noindent {\bf Theorem}(\cite{KS1} 2.6.4.)  {\it Let $M$ be a complete Riemannian
manifold with finite volume, and $X$ be an NPC space.  Let $\Gamma = \pi_1(M)$,
and $\rho : \Gamma \rightarrow {\rm Isom}(X)$ be a homomorphism. Suppose that 
there exists  a $\rho$-equivariant map $\phi: \tilde{M} \rightarrow \Tbar$ with 
finite energy $E$. Then there exists an $\rho$-equivariant minimizing sequence $\{ u_i \}
: \tilde{M} \rightarrow \Tbar$ so that for any compact set $K \subset M$ and $i$ sufficiently
large (depending on $K$) the $u_i$'s are Lipschitz continuous on (the lift to $\tilde{M}$
of) K, with the pointwise Lipschitz constant bounded by $C(x)E^{1/2}$, where $C(x)$ is a
locally bounded function depending only on $M$. }\\

This gives us a minimizing subsequence whose energy is converging to $E_0 = \inf_{u} \{ E(u): 
 \mbox{ $u$ $\rho$-equivariant} \}$. To ensure that the images of the sequence are not escaping toward infinity (the geometric boundary of $X$), we need to find a way to anchor the images
within the finite region. The following statement~\cite{KS1} 2.6.5. provides such a criterion.

\noindent {\bf Theorem}(~\cite{KS1} 2.6.5.) {\it Let $M$ be a complete Riemannian
manifold, possibly with compact Lipschitz boundary $\partial M$.  Let $\Gamma = \pi_1(M)$,
and $\rho : \Gamma \rightarrow {\rm Isom}(X)$ be a homomorphism.  Let $u_i$ be
an equivariant minimizing sequence which has local modulus of continuity control.  That is,
for each $x \in \tilde{M}$ we assume there is an equivariant function $w(x, r)$,  $ 0 \leq r
< r_x$ which is monotone increasing in $r$, which satisfies $w(x, 0) = 0$ and 
\[
\sup_{i} \sup_{ d(x, z) \leq r} d(u_i(x), u_i(z)) \leq w(x, r).
\]
Then the sequence converges locally uniformly to an equivariant harmonic map $u$ if and only
if there exists an $x \in \tilde{M}$ at which the sequence $u_i(x)$ is convergent.
}\\

Let $u_i$ be the energy minimizing sequence whose existence is given by the first theorem quoted.
Then recall we have the uniform Lipschitz bound for $\{ u_i \}$, which in turn ensures the 
existence of $w(x, r)$ in the statement of the second quoted theorem.  We now claim that
there is a subsequence, which we will denote by  $\{ u_i \}$ again, such that for some $x$,
$\{ u_i (x) \}$ is convergent.  

Recall the representation  $\rho : \Gamma \rightarrow \map$  induces a proper action on
$\Tbar$.  By the same argument as in the proof of the existence theorem for the compact 
domain, we know that there exists some $C >0$ such that 
\[
d_{wp} (\gamma_{\alpha} u_i(x), u_i(x)) < C 
\]
for $x \in F$ (recall $F$ is a fundamental region in $\tilde{M}$), and for
all $i$.  In particular we know $u_i [F]$ is contained in a sublevel set
of the functional $\delta$, which is bounded by the properness of $\rho$.
Hence for any $x \in F$, there exists some convergent subsequence $\{ u_i(x) \}$.

\end{pf}

\section{Proof of the Uniqueness Theorem}

\begin{thm}[Uniqueness of Harmonic Maps]
The harmonic map  is unique
within the class of finite energy maps which are $\rho$-equivariant, provided
that the map is neither constant nor into a geodesic in $\Tbar$.
\end{thm}

\begin{pf}
It was shown in~\cite{KS1} (2.2v) that if $u$ and $v$ are two points in $Y= L^2_{\rho} (M, X)$
where $X$ is an NPC space, and if $u_t$ is a interpolation between $u$ and $v$ given
by the geodesic homotopy, that is, $u_t(x)$ is a point on the unique geodesic connecting
$u(x)$ and $v(x)$ with $t \in [0, 1]$ parametrizing the geodesic with a constant speed,
then
\begin{eqnarray}
d^2(u_t(x), u_t(y)) & \leq & (1-t)d^2(u(x), u(y)) + t d^2(v(x), v(y))\\
 &  & - t(1-t)[d(u(y), v(y)) - d(u(x), v(x))]^2. 
\end{eqnarray}
We first make the statement infinitesimal (by making $x$ and $y$ close), in order to define
the energy densities $| \nabla u |^2$ and $| \nabla v |^2$ weakly.  Then integrate 
over $M$ the inequality
among the 
energy densities to obtain the following inequality among 
energy functionals $E(u), E(v)$ and $E(u_t)$ 
\[
E(u_t) \leq (1-t) E(u) + t E(v) - t(1-t) \int_M | \nabla d(u,v) |^2.
\]
Clearly this says that $E(u_t)$ is a  convex function in $t$.  

Suppose that $u, v : \tilde M \rightarrow \Tbar$ are both energy minimizing maps 
which are $\rho$-equivariant.  Recall the definition 
$E_0 = \inf \{ E(u): \mbox{$u$\ $\rho$-equivariant} \}$.  Then by the
construction of the energy minimizing maps $u$ and $v$, we know 
\[
E(u) = E(v) = E_0. 
\]
By the convexity along the geodesic homotopy connecting $u$ and $v$, this
implies that $\int_M | \nabla d(u,v) |^2 = 0$, since the function $E(u_t)$
cannot be strongly convex in this situation.  This in turn implies that 
the two maps are {\it parallel}, that is, $d(u(x), v(x)) \equiv C$ where
$C$ is a constant.

Now  define a map $\Phi (u, v) : \tilde{M} \rightarrow \Tbar \times \Tbar$ 
so that
$\Phi (u, v) (x) = (u(x), v(x))$. Then since both $u$ and $v$ are (Lipschitz) 
continuous, so is $\Phi (u, v)$.

 We now look at the inverse image of 
the diagonal $D = \{ (p, q):
p=q \}$ of $\Tbar \times \Tbar$.  Since $D$ is closed in $\Tbar \times \Tbar$,
the inverse image of it by $\Phi(u,v)$ is closed in $M$.  There are two cases
to consider; the first is when $[\Phi(u, v)]^{-1}(D)$ is $M$, in which case
$d(u(x), v(x)) \equiv 0$, or equivalently $u \equiv v$, the second case
being $M \setminus [\Phi(u, v)]^{-1}([\Tbar \times \Tbar] \setminus D)$ is 
nonempty open set of $M$. Note 
that the first case implies the uniqueness.

We proceed to analyze the second case.  We will show that if $u$ and $v$ do
not coincide, then there is a nonempty open set $U$ in $M$ which is mapped
to $[{\T}_C \times {\T}_{C'}] \setminus D$ for some subsets $C$ and $C'$ of $\cal S$
where each represents a collection of $N$ mutually disjoint simple closed 
geodesics on the surface $\Sigma$ with $0 \leq N < 3g-3$.

We first look at the inverse image $[\Phi(u, v)]^{-1}([\T \times \T] \setminus D)$.  Since  $[\T \times \T] \setminus D$ is open in $[\Tbar \times \Tbar]$,
the inverse image is an open set in $M$.  If the set is nonempty,
then the inverse image is the open set $U$ we had set to find.
If it is empty (that is, if the image of $\Phi (u, v)$ does not intersect
with $[\T \times \T] \setminus D$), then we proceed to look at the inverse 
image of 
$[{\T}_c \times {\T}_{c'}] \setminus D$ by $\Phi (u, v)$ with $c$ and $c'$
being  closed geodesics on $\Sigma$, i.e.  subsets of $\cal S$ whose 
cardinality is one (which we will denote by $|c|= |c'|=1$.)  
Since the set
 $[{\T}_c \times {\T}_{c'}] \setminus D$ is an open set in $[(\Tbar \setminus
\T) \times (\Tbar \setminus \T)] \backslash D$ with the relative topology
induced by the topology of $\Tbar$, it follows that the inverse image
$[\Phi(u, v)]^{-1} \Big([{\T}_c \times {\T}_{c'}] \setminus D \Big)$ is open in $M$.
If this set is nonempty for some $c$ and $c'$ in $\cal S$, then we have the 
open set as desired.  If not, (that is, if the image of $\Phi (u, v)$ 
does not intersect with $\Big[ \Big( {\T} \cup \{ \cup_{c \in {\cal S}} {\T}_c 
\} \Big) \times
 \Big( {\T} \cup \{ \cup_{c' \in {\cal S}} {\T}_{c'} \} \Big) \Big] \setminus 
D$, then 
we iterate the previous argument, this time applied to the 
open subset $[{\T}_C \times {\T}_{C'}] \backslash D$ of $\Big[ \Big(
\Tbar \setminus
\{ {\T} \cup \{ \cup_{c \in {\cal S} } {\T}_c \} \Big) \times \Big( 
\Tbar \setminus
\{ {\T} \cup \{\cup_{c \in {\cal S}} {\T}_c \} \Big) \Big] \setminus D$ with the relative 
topology. Here $C, C'$ are taken to be subsets of $\cal S$ of cardinality
two (representing pairs of two mutually disjoint closed geodesics.)

Recall that the cardinality $|C|$ of a subset $C$ of $\cal S$ representing 
a set of mutually
disjoint simple closed geodesics has an upper bound of $3g - 3$.
Hence the iterative process as described above has to stop sometime.
Also note that when $3g-3$ geodesics pinches off, then the original
surfaces has been disconnected into a collection of pairs of
pants, (or equivalently thrice punctured spheres), which is to say
${\T}_C$ is a point in $\Tbar$ if $|C|$ is maximal
of $3g-3$. This is due to the fact that there is no deformation
of the conformal structure of the thrice punctured sphere.  
This implies that both $u$ and $v$ are constant maps
into $\Tbar$, which violates the hypothesis that the maps are
of rank bigger than one at some point.

Therefore the  cases when the uniqueness may not hold
is when there is a nonempty  proper open set $U$ of $M$ which is 
mapped by $\Phi (u, v)$
to $[{\T}_C \times {\T}_{C'}] \setminus D$, where $C$ and $C'$ are 
two subsets of $\cal S$ of the same cardinality $|C| = |C'| = N$
with $0 \leq N < 3g-3$, each 
representing a set of $N$ mutually disjoint simple closed geodesics.

We know that the map $u_t$ defined by interpolating $u$ and $v$
via geodesic homotopy is also energy minimizing, for $E(u_t) = E_0$, $t \in [0, 1]$.

It is proven in~\cite{Y2} that any open \weil\ geodesic segment lies
in one single component ${\T}_C$ of $\Tbar$.  In our situation, 
we have $u(x)$ in ${\T}_C$ and $v(x)$ in ${\T}_{C'}$ when $x$ is in $U$.
Hence the open geodesic segment with its end-points $u(x)$ $v(x)$
lies in a component ${\T}_{\tilde{C}}$.  It is also shown in~\cite{Y2} 
that if $|C| = |C'|$ and $C \neq C'$, then ${\T}_C$ and ${\T}_{C'}$
are either disjoint or intersect normally.
If $C = C'$ then ${\tilde{C}} = C = C'$, otherwise 
$| \tilde{C} | > |C| = |C'|$.  This tells us that $u_{1/4}(x)$
and $u_{3/4}(x)$ are both in ${\T}_{\tilde{C}}$ for each $x$ in $U$.
Denote $u_{1/4}$ by $f$ and $u_{3/4}$ by $g$.

Recall~\cite{KS1} that an energy minimizing map $u$ is locally the solution of 
the Dirichlet problem, i.e. given an open domain $\Omega$ with a Lipschitz 
boundary in $M$ in $U$, the restriction
of the map $u$ to $\Omega$ is the solution to the Dirichlet problem whose
Dirichlet data is given by $u \big|_{\partial \Omega}: \partial \Omega \rightarrow 
\Tbar$.  
Moreover when the target is smooth, then the map $u$ is smooth on 
$\Omega \subset U$.   
Hence in our setting where the target of $f$ and $g$ is the smooth manifold
${\T}_{C'}$, we know that $f$ and $g$ are both smooth on
any open set properly contained in $U$.  

We now need the following calculation as demonstrated in~\cite{SY}.

Let $\rho(x) = d(f(x), g(x))$ defined on $U \subset M$. Let $1 \leq \alpha 
\leq  m$
be the induces for an orthonormal frame $\{e_{\alpha} \}$ at $x$ in $U$.
And let $1 \leq i \leq n$ be the induces for  orthonormal frames 
on ${\T}_{\tilde{C}}$ at $f(x)$ and $g(x)$.  Also denote by subscripting $1 \leq  \overline{i} \leq n$
 differentiation with respect to the last $n$ coordinates of the 
\weil\ distance function
$d: \Tbar \times \Tbar \rightarrow {\bf R}$.   
  Then we have
\[
\triangle \rho^2 = 2 \Sigma_{\alpha, i} (d_i (f_{i})_{\alpha} + 
d_i (g_{i})_{\alpha})^2 + 2 \Sigma_{\alpha} \rho d_{X_{\alpha} X_{\alpha}}
\]
where $d_{X_{\alpha} X_{\alpha}}$ is the Hessian of the distance function $d$.
Since $\rho$ is constant on $U$, the equation above implies that
\[
\Sigma_{\alpha} \rho d_{X_{\alpha}X_{\alpha}} \equiv 0.
\]
where the vector field $X_{\alpha}$ on $T_{f, g} \big( {\T}_{\tilde{C}} \times 
{\T}_{\tilde{C}} \big)$ is given by $(f_{*}e_{\alpha}, g_{*}e_{\alpha})$.

While we have the following second variation formula of the pulled-back
distance function;
\[
\Sigma_{\alpha} d_{X_{\alpha} X_{\alpha}} =  
\Sigma_{\alpha} \int_0^1 \{ \| \nabla^{\perp} V_{\alpha} \|^2 - 
< R^{wp} (V_{\alpha}, T)T, V_{\alpha}> \} dt 
\]
where $V_{\alpha}$ is the Jacobi field along the geodesic
segment $\overline{f(x)g(x)}$ (here parametrized by $t \in [0, 1]$
whose value at $f(x)$ is $f_{*}e_{\alpha}$ and whose value at
$g(x)$ is $g_{*}e_{\alpha}$.  $V^{\perp}$ denotes the component 
of $V$ normal to the geodesic $\overline{f(x)g(x)}$. And $T$ 
denotes the tangent vectors along the arc-length
parametrized (by $t \in [0, 1]$) geodesic $\overline{f(x)g(x)}$.

It has been known that the sectional curvature of \weil\
metric is strictly negative (see \cite{Wo2} for example.)  Thus it follows that
\[
\Sigma_{\alpha} d_{X_{\alpha} X_{\alpha}} > 0,   
\]
unless the image if $u$ and $v$ are contained in the 
geodesic.  This contradicts with the vanishing of $\triangle \rho^2$.  Hence
we have to have $u \equiv v$ on $U$, which then imposes that
\[
d(u(x), v(x)) \equiv 0
\]
on $M$.

\end{pf}

\section{Proof of the Stratification Property}

We are now ready to prove the stratification theorem of harmonic images.

\begin{thm}[Stratification Property of Harmonic Images]
Suppose that $\rho$ induces a proper action on $\Tbar$ and that there
exists a unique energy minimizing harmonic map $u: \tilde{M} \rightarrow \Tbar$.
The subgroup $\Gamma = \rho (\pi_1 (M))$ is reducible by a collection $C$
of mutually disjoint closed geodesics, if and only if  the image of the energy minimizing map
$u$ lies entirely in the boundary component $\overline{{\T}_C}$. 
\end{thm}

\begin{pf}
It is clear that if the image lies on some $\overline{{\T}_C}$, then the map $u$ 
is equivariant with respect to a representation which is reducible by $C$.

We consider the other direction; suppose that the representation is reducible.

Recall that the boundary ${\T}_C$ with $C= \{ c_i \}_{i=1, ...,n} \subset {\cal S}$
where $c_i$'s are mutually disjoint closed homotopically nontrivial simple closed curves,
is the \teich\ space of the nodal surface ${\Sigma}_C$ obtained by pinching all of the $c_i$'s.
If the nodal surface consists of more than one path-connected components, which
we will denote by ${\Sigma}_{C_{\alpha}}$ with $\cup_{\alpha} C_{\alpha} = C$, ${\T}_C$
is a product space of the \teich\ spaces ${\T}_{\alpha}$ of different components 
${\Sigma}_{C_{\alpha}}$.

In the case of ${\T}_C$ being a direct product of ${\T}_{\alpha}$'s with no two
components ${\Sigma}_{C_{\alpha}}$'s are diffeomorphic, then the restriction
of the mapping class group $\map$ is the direct product of the mapping class groups
${\rm Map}({\Sigma}_{C_{\alpha}})$ of the path-connected components ${\Sigma}_{C_{\alpha}}$'s.

In the case when there is at least a pair of path-connected components $ {\Sigma}_{C_{\alpha}}$
which are diffeomorphic to each other, then the restriction of $\map$
contains a subgroup of a symmetric group of $N$ elements, where $N$ is the 
number of the path connected components of ${\Sigma}_C$.

Note that the completion ${\Tbar}_C$ of ${\T}_C$ is an NPC space {\it on its own}, 
equipped with the 
\weil\ distance function, whose
isometry group is the restriction of the mapping class group $\map$ on ${\Sigma}_C$.

Since $\rho$ induces a proper action on $\T$, the sublevel sets of the functional
$\delta: \T \rightarrow \bf R$ is bounded in $\T$, and hence is bounded
in ${\T}_C$.  Since the $\rho [\pi_1 (M)] \subset \map$ is a reducible subgroup, 
it is a subgroup of the isometry group of ${\T}_C$, which is the restriction of
$\map$ on ${\T}_C$.  Now we have a reduced setting as follows.

Let $\rho': \pi_i (M) \rightarrow {\rm Map}({\Sigma}_C)$ induced by the restriction
of $\rho : \pi_1 (M) \rightarrow \map$. We know $\rho'$ is a  proper action
on ${\Tbar}_C$, and by the 
Existence Theorem 1, we know that there exists a $\rho'$-equivariant harmonic map
$u' : \tilde{M} \rightarrow {\Tbar}_C$.  Now by the uniqueness,
there is no other $\rho'$-equivariant harmonic map, therefore no other 
$\rho$-equivariant map. 

\end{pf}

\section{Unique Continuation Property of Harmonic Images}

\begin{thm} Suppose that a Riemannian region $\Omega$ is mapped into ${\Tbar}$ 
by a geometrically finite energy-minimizing map $u$, and that there is a point 
$p$ in $\Omega$ mapped into a \teich\ space ${\T}_C \subset {\T}$ for some collection of
mutually disjoint simple closed curves.  Then  the entire
image $u[\Omega]$ lies in ${\T}_C$.
\end{thm}

\begin{pf}[of the Theorem] Suppose the contrary, that is, there exists a point $p$ in $\Omega$ 
which is mapped to a point $u(p)$ in $\overline{{\T}_c}$ where $c$ is a homotopically
nontrivial simple closed curve.  Let $u_1$ be
the coordinate function near $u(p)$ such that $u_1$ vanishes
on $\overline{{\T}_c}$, as defined in~\cite{Y2}.  We will show then that the entire image
of $\Omega$ by $u$ lies in $\overline{{\T}_c}$, a contradiction to
a hypothesis of the statement of the theorem.

The harmonic map equation 
for $u : \widetilde{S^2 \backslash \{p_i\}}  \rightarrow \Tbar$
for the first coordinate function $u^1$ is
\[
\triangle_{\Omega} u^1 + \Gamma^1_{\alpha \beta}(u(x)) u^{\alpha}_i u^{\beta}_j 
g^{ij} (x) = 0
\]
where $\Gamma^1_{\alpha \beta}$ is the Christoffel symbol for the
\weil\ metric evaluated at $u(x)$.  Note that $u^1 (x)$ is nothing but the pulled-back
function $(u^* u_1)(x)$.  (Also note $u^2 (x) = (u^* \theta_1)(x)$ where $\theta_1$
is the Fenchel-Nielsen angle coordinate around $c$.) We will show that
near $u(p)$ in $\overline{{\T}_c}$, $u^1$ satisfies the
differential inequality weakly;
\[
\triangle_{\Omega} u^1 \leq C u^1
\]
for some $C >0$.
We do so by using the harmonic map equation for $u^1$ and by 
showing that the nonlinear terms satisfy
\[
\sum_{\alpha, \beta} \Gamma^1_{\alpha \beta}(u(x)) u^{\alpha}_i u^{\beta}_j 
g^{ij} (x) < C u^1 (x)
\] 
for some $C > 0$. In fact, we will show 
\[
 \sum_{\alpha, \beta} \Gamma^1_{\alpha \beta}(u(x)) u^{\alpha}_i u^{\beta}_j 
 g^{ij} (x) 
= O (u^1 (x))
\]
which allow us to choose some $C > 0$ for sufficiently small $u^1$.

We will use the following result from~\cite{Y2} that with respect to the coordinate 
system $t = (u_1, \theta_1, t_2,...,t_{3g-3})$
near $t=0$, the \weil\ metric tensor has the following form;
\[
\begin{array}{lll}
G_{11}(t) & = & 1 + O( (u_1)^4 ) \\
G_{12}(t) & = & 0 \\
G_{1j}(t) & = & O( (u_1)^3 ) \ \ (j > 2) \\
G_{22}(t) & = & \frac{1 + O( (u_1)^4 )}{4} (u_1)^6 \\
G_{2j}(t) & = & O((u_1)^6) \ \ (j > 2) \\
G_{ij}(t) & = & (1 + O((u_1)^4)) G_{ij}(0, 0, t_2,...,t_{3g-3}) \ \ 
\mbox{ for $i, j > 2$}.
\end{array}
\]
The Christoffel symbols are obtained from the metric tensor by the
following formula (see~\cite{GDH} for example)
\[
\Gamma^1_{ij} = \frac{1}{2} \sum_{l} G^{1l} (G_{il, j} + G_{lj, i} - G_{ij,l}).
\]
Then we have 
\[
G^{11} = O(1), G^{12} = O(1) \mbox{ and } G^{1j} = O((u_1)^3) \ \ 
\mbox{for $j >2$}
\]
and consequently we have
\[
\Gamma^1_{11} = O((u_1)^3), \ \ \Gamma^1_{12} = O((u_1)^4), \ \ \Gamma^1_{22}
= O((u_1)^5)
\]
\[ 
\Gamma^1_{1j} = O((u_1)^2), \ \ \Gamma^1_{2j} = O((u_1)^3) \ \
\mbox{ and } \ \ \Gamma^1_{ij} = O((u_1)^3) \mbox{ for $i, j > 2$}.
\]
It was shown in~\cite{GS} (Theorem 2.4) that the energy density of the 
energy minimizing map is bounded by the total energy of the map, and 
hence that the map is uniformly Lipschitz continuous.  This implies that
$u^{\alpha}_i$ is bounded unless $\alpha = 2$, for the \weil\ metric
tensor degenerates in the direction of $\del / \del \theta$ as $u_1$
goes to zero.  However
we know that $u$ is geometrically finite and hence $u^2 = u^* \theta = O (1)$
as $u_1$ vanishes, which in turn implies that  $u^2_i = O ( (u_1)^{-1} )$.
Note then that 
\[
\Gamma^1_{2 \beta} u^2_i u^2_j = O( (u_1)^3 ) O( (u_1)^{-1} ) O(
 (u_1)^{-1} ) = O( u_1 ) = O( u^1(x) )
\]
which is what we needed to show the inequality 
$\triangle u^1 \leq C u^1.$
for $u^1 > 0$.  

Note that in the region $u^1 > 0$, the energy minimizing 
map is smooth, for $\T$ is an open simply-connected negatively-curved
manifold.  To see that this inequality holds distributionally
for the region $u^1 \geq 0$, that is, for a compactly supported 
non-negative test function $\phi$ on $\Omega$
\[
\int_{\Omega} (u^1) \triangle \phi d \mu \leq \int_{\Omega} 
C (u^1) \phi d \mu
\]
we integrate by parts.  Define ${\Omega}_1 = \{ x \in \Omega : u^1 (x) > 0 \}$
and ${\Omega}_2 = \{ x \in \Omega: u^1 (x) = 0 \}$
\[
\begin{array}{lll}
\int_{\Omega} (u^1) \triangle \phi d \mu & = & \int_{\Omega_1}
(u^1) \triangle \phi d \mu + \int_{\Omega_2}
(u^1) \triangle \phi d \mu  \\
 & = & - \int_{\Omega_1} \langle \nabla (u^1), \nabla \phi \rangle d \mu 
+ \ \int_{\del \Omega_1} (u^1) \frac{\del \phi}{\del \eta} d \sigma \\
 & = & \int_{\Omega_1} \triangle (u^1) \phi d \mu - \int_{\del \Omega_1} \phi 
 \frac{\del u^1 }{\del \eta} d \sigma \\
 & \leq & \int_{\Omega_1} C ( u^1 ) \phi d \mu  
\end{array}
\]
where $\eta$ is the outer unit normal vector to the boundary, and 
$d \sigma$ is the induced measure on the boundary.  The second equality
is due to an integration by parts.  Note here that $u^1 \equiv 0$ 
on the boundary $\del \Omega_1$.  The third equality is due to 
an integration by parts and a claim $\frac{\del u^1 }{\del \eta} \equiv 0$
on $\del \Omega_1$.  To see that the claim holds, suppose that $x$ is 
a point on $\del \Omega_1$ and that $\frac{\del u^1 }{\del \eta}(x) \neq 0$. 
Recall from~\cite{GS} 
the construction
of the homogeneous approximating map $u_*$ from a tangent space $T_x \Omega$
into the tangent cone $C_{u(x)} \Tbar$, which is obtained by taking 
a limit of $u_{\lambda, \mu} {\mu}^{-1} u (\lambda x) : B_{{\lambda}^{-1} 
\Omega} \rightarrow {\mu}^{-1} \Tbar$ where ${\mu}^{-1} \Tbar$ is the 
NPC space whose distance measured from $u(x)$ is multiplied by ${\mu}^{-1}$,
and where $B_{{\lambda}^{-1} \Omega} (x)$ has a scaled metric 
$g_{\lambda}(y) = g (\lambda y)$  ($\lambda y$ means the dilation
on the tangent space which is identified with $\Omega$ via the exponentially map based at $x$.)  We let $\lambda$ and $\mu$ tend to zero as described 
in~\cite{GS} to get a limiting map $u_*$, which we call an homogeneous
approximating map.  It is itself an energy-minimizing map 
from ${\bf R}^n = T_x \Omega$ to the tangent cone $C_{u(x)}\Tbar = {\bf R}^+  \times {\bf R}^{(6g-6-2)} 
$. Note that the Fenchel-Nielsen twist coordinate $\theta_1$ has dropped
out when one takes the limits of dilations.  We will still denote the coordinates of 
${\bf R}^+ \times {\bf R}^{(6g-6-2)}$ by $(u_1, x_3,..., x_{6g-6})$.
Then $u_1$ is a linear coordinate function with respect to
the distance function induced by the scaling, and therefore the pulled-back function 
$(u_*)^* u_1$ on ${\bf R}^n$ is superharmonic, thus satisfying the
mean-value inequality
\[
\inf_{B_r (x)} (u_*)^* u_1 \leq \frac{1}{r^n} 
\int_{B_r (x)} (u_*)^* u_1 d \mu.
\]
Since $x \in \Omega$ is mapped to $\overline{{\T}_c} = \{ u_1 = 0 \}$.
the approximating homogeneous map $u_*$ at $x$ maps the origin of
${\bf R}^n$ to the origin of ${\bf R}^+  \times {\bf R}^{(6g-6-2)}$. Then
the mean-value inequality says that the entire image of $u_*$ has to
lie in the set $\{ u_1 = 0 \} \cong {\bf R}^{(6g-6-2)}$.   Now note that this
contradicts to the condition $\frac{\del u^1 }{\del \eta}(x) \neq 0$, hence the
claim is shown.

Therefore we have the desired inequality 
\[
\triangle u^1 \leq C u^1.
\]
satisfied distributionally on a Riemannian domain $\Omega$.   The 
DiGiorgi-Nash-Moser estimate (see, for example, ~\cite{GT}) then implies that
\[
\inf_{U} u^1 \geq \tilde{C} \int_\Omega u^1 d \mu
\]
for any proper open subset $U$ of $\Omega$, where $\tilde{C} > 0$ depends on
$U$ and $\Omega$ but independent of the function.  

Now suppose that a point $x$ is mapped into ${\T}_c$, thus we have
$u^1(x) = 0 = \inf_U u^1$ and the last
inequality above then implies that $u^1 \equiv 0$ on the entire domain 
$\Omega$.  

This concludes the proof of the
theorem.

\end{pf}

\end{document}